\documentclass[a4paper]{article}
\usepackage{color}
\usepackage{caption}
\usepackage{float}
\usepackage{amsmath} 
\usepackage{calc}
\newlength{\depthofsumsign}
\usepackage{amsfonts,amsmath,amssymb,amsthm}
\usepackage{graphicx}
\usepackage{mathrsfs}
\usepackage[title, titletoc, header]{appendix}

\newtheorem{definition}{Definition}[section]

\newtheorem{theorem}{Theorem}[section]

\newtheorem{problem}{Problem}[section]

\usepackage{amsmath}

\newtheorem{remark}{Remark}[section]

\begin{document}

	\begin{center}
		\large{ \textbf{Optimal selection of local approximants in RBF-PU interpolation}}
	\end{center}

	\begin{center}
		Roberto Cavoretto$^a$, Alessandra De Rossi$^{a}$, Emma Perracchione$^{b}$ 
	\end{center}
	
	\begin{center}
		$^a$ Department of Mathematics \lq\lq G. Peano\rq\rq, University of Torino, via Carlo Alberto 10, I--10123 Torino, Italy		\\
		$^b$ Department of Mathematics \lq\lq T. Levi-Civita\rq\rq, University of Padova, via Trieste 23, I--35121 Padova, Italy
	\end{center}
	\vskip 0.5cm

	\textbf{Abstract.}
				The Partition of Unity (PU) method, performed with local Radial Basis Function (RBF) approximants, has been proved to be an effective tool for solving large scattered data interpolation problems. However, in order to achieve a good accuracy, the question about how many points we have to consider on each local subdomain, i.e. how large can be the local data sets,  needs to be answered. Moreover, it is well-known that also the shape parameter affects the accuracy of the local RBF approximants and, as a consequence, of the PU interpolant. Thus here, both the shape parameter used to fit the local problems and the size of the associated linear systems are supposed to vary among the subdomains. They are selected by minimizing an \emph{a priori} error estimate. As evident from extensive numerical experiments and applications provided in the paper, the proposed method turns out to be extremely accurate also when data with non-homogeneous density are considered.\\

\section{Introduction}
\label{Intro}

Given a set of multivariate data, we aim at finding a function that accurately fits such points. This problem is rather common in many applied sciences, such as in physics, biology, geophysics and Earth's topography. Moreover, dealing with applications, one often faces the problem of approximating large and \emph{irregular} data sets, i.e. data which are far from being uniform or quasi-uniform. In these cases, since problems as lack of information, i.e.  low density of points, or ill-conditioning, i.e. high density of data, arise, the fitting process becomes a challenging computational issue.

Because of the above mentioned problems, recently, the approximation theory has driven its attention on local techniques. Specifically, the approximation of irregularly distributed data via local schemes has gained much attention in both meshfree and mesh-dependent methods. For the latter, the problem results particularly hard and the choice of the mesh turns out to be crucial. As example, in \cite{Davydov1}, in order to build the local approximating fits, the authors consider spline functions on a uniform triangulation with $C^1$ or $C^2$ continuity. While, in \cite{Bozzini}, polyharmonic splines are effectively used to fit irregular and truly large data sets.

Another major class of techniques, which are known in literature as meshfree or meshless methods, includes RBF approximants \cite{Buhmann03,Iske11}. They obviously take advantage of being independent from the mesh and, as a consequence, they are easy to implement in any dimension. Indeed, a local hybrid approximation technique for data with non-homogeneous density, obtained by means of both splines and RBFs, is presented in \cite{Davydov2}. But, since bivariate spline functions are involved, the method again strongly depends on the mesh. To avoid this drawback, we focus on purely meshless methods. In this context, the scattered data problem of huge and irregular sets of points is usually performed by means of least squares approximation \cite{Schaback00a,Wendland01}. Here instead, our fitting criterion consists in exactly matching the measurements at their corresponding locations, i.e. we focus on interpolation. To this aim, the PU method  performed with local RBF interpolants turns out to be really meaningful \cite{Babuska97,Fasshauer}.

The basic idea of the PU technique consists in decomposing the domain into several \emph{subdomains} or \emph{patches} forming a covering of the original domain. When the PU method is applied in the context of interpolation, except for particular cases \cite{Safdari}, such subdomains  are always supposed to be hyperspheres of a fixed size \cite{Cavoretto15b,Fasshauer,Heryudono16,Shcherbakov}. But, in case of irregular data this might lead to inaccurate approximations. In \cite{Safdari} variable subdomains are used for an \emph{ad hoc} 2D problem in finance; specifically, even if data are not quasi-uniform, they have a precise and well-known structure. Thus, the PU subdomains are constructed following exactly their distribution. 

Our aim consists instead in developing a method which enables us to select, independently from the node distribution, suitable sizes of the different PU subdomains. Furthermore, we also take into account the critical choice of the shape parameter of the basis function. In fact, it can greatly influence the accuracy of final fit. 

To such scope, we compute subsequent a priori error estimates depending on both  the shape parameter and  the size of the PU subdomain. Then, for each patch we select  the \emph{optimal} couple of values, i.e. the subdomain size and the shape parameter, used to solve the local interpolation problem. The error estimates are found out via a \emph{modified} Leave One Out Cross Validation (LOOCV) scheme \cite{Fasshauer15,Golub79,Rippa}. More precisely,  since our problem depends on two quantities, for each patch  we  perform a Bivariate LOOCV (BLOOCV).
The resulting method, named BLOOCV-PU, turns out to be extremely accurate compared with the classical PU technique. This improvement, in terms of accuracy, becomes particularly meaningful when data with non-homogeneous density are considered.

The complexity of the algorithm is also taken into account. Specifically, the computational issue consisting in finding all the points belonging to a given subdomain is performed with the use of a novel data structure, the so-called Integer-based Partitioning Structure (I-PS). It leads to a considerable saving in terms of computational time with respect to the most advanced searching techniques \cite{Cavoretto15b,AIP_15Emma}.

Numerical experiments show the good performances of the BLOOCV-PU meth-\ od in case of quasi-uniform data, and underline the benefits of such a flexible approach with irregular points. Moreover, we investigate two applications with real world data, including a  benchmark glacier
data set and points with highly varying densities describing a terrain.

The guidelines of the paper are as follows. In Section \ref{PUM_th}, after briefly reviewing the main theoretical features of the PU method, we  introduce the BLOOCV-PU interpolant. The computational aspects of such algorithm and its complexity are described in Sections \ref{comput_asp} and \ref{analisi_costi}, respectively. Extensive numerical experiments and applications with real world data, carried out in Sections \ref{Numerica_ex} and \ref{foresta} respectively, are devoted to test the accuracy and the flexibility of the BLOOCV-PU approximant. Finally, in Section \ref{rem_comcl}, we deal with conclusions and work in progress.

We point out that the  \textsc{Matlab} software is made available to the scientific community in a downloadable free  package: 
\begin{center}
	{\tt http://hdl.handle.net/2318/1559094}.
\end{center}

\section{Formulation of the BLOOCV-PU interpolant}
\label{PUM_th}
In Subsection \ref{PUM} we first briefly review the main theoretical aspects concerning the PU interpolation and then in Subsection \ref{opt_ep} we focus on the local selection of suitable sizes of the patches and shape parameters. In order to achieve such scope, for each local interpolation problem we need to determine an error estimate depending on both  the size of the subdomain and the shape parameter of the basis function. 
\subsection{The partition of unity method}
\label{PUM}

The scattered data interpolation problem consists in recovering a function $f: \Omega\to \mathbb{R}$, $\Omega$ being a bounded set in $\mathbb{R}^M$, using a set of samples of $f$ on $N$ distinct data points or nodes ${\cal X}_N= \{\boldsymbol{x}_i, i=1,\ldots,N\}\subset\Omega$, namely $\boldsymbol{f} =(f_1, \ldots , f_N)^T$, $f_i = f(\boldsymbol{x}_i)$, with $\boldsymbol{x}_i\in {\cal X}_N$. More precisely, we aim at solving the above mentioned problem for truly large data sets. To this scope, the PU method, first introduced in \cite{Babuska97,Melenk96}, is a well-suited numerical tool.

The idea behind the PU method is to start with a partition of the open and bounded domain $ \Omega$ into $d$ subdomains $ \Omega_j$,  such that $ \Omega  \subseteq \cup_{j=1}^{d} \Omega_j$,  with some mild overlap among them \cite{Babuska97,Fasshauer,Heryudono16,Melenk96,Wendland02a}. 
Together with these subdomains, we need  a family of compactly supported, non-negative, continuous functions $W_j$, $j=1,\ldots,d$, which form a $k$-stable partition of unity, i.e.
\begin{equation*}
\sum_{j=1}^{d} W_j(\boldsymbol{x}) = 1, \quad \boldsymbol{x} \in \Omega,
\end{equation*}
and for every $ \boldsymbol{\beta} \in \mathbb{N}^{M}$, with $| \boldsymbol{\beta} |	 \leq k$, there exists a constant $C_{ \boldsymbol{\beta}} >0$ such that
\begin{displaymath}
\left\| D^{ \boldsymbol{\beta}} W_j \right\|_{L^{ \infty} ( \Omega_j)} \leq  \dfrac{C_{ \boldsymbol{\beta}}}{ \left(  \sup_{ \boldsymbol{x}, \boldsymbol{y} \in \Omega_j} \left\| \boldsymbol{x}-\boldsymbol{y} \right\|_2 \right)^{ | \boldsymbol{\beta}|}}, \quad j=1, \ldots, d. 
\end{displaymath}

More precisely, in what follows  we consider the so-called Shepard's weights which are defined as
\begin{equation*}
W_j(\boldsymbol x) = \frac{ \displaystyle \bar{W}_j (\boldsymbol x)}{  \displaystyle \sum_{k=1}^{d} \bar{W}_k (\boldsymbol x)}, \quad j=1, \ldots,d,
\end{equation*}
where $\bar{W}_j $ are compactly supported functions forming a partition of unity. 

Once we choose the partition of unity $ \{ W_j \}_{j=1}^{d}$, the global interpolant is formed by the weighted sum of $d$ local approximants $R_j$, i.e.
\begin{equation}
{\cal I}\left( \boldsymbol{x}\right)= \sum_{j=1}^{d} R_j\left( \boldsymbol{x} \right) W_j \left( \boldsymbol{x}\right), \quad \boldsymbol{x} \in \Omega.
\label{intg}
\end{equation}

In particular here $R_j$ denotes a  RBF interpolant defined on a subdomain $\Omega_j$ of the form
\begin{eqnarray}
\label{rad1}
R_j( \boldsymbol{x} )= \sum_{k=1}^{N_j} c_k^j \phi ( || \boldsymbol{x} -  \boldsymbol{x}^j_k ||_2 ), 
\end{eqnarray}
where $\phi:[0,\infty) \rightarrow \mathbb{R}$ is called RBF, $||\cdot||_2$ denotes the Euclidean norm,  $N_j$ indicates the number of data points belonging to $\Omega_j$ and $\boldsymbol{x}_k^j \in {\cal X}_{N_j}= {\cal X}_N \cap \Omega_j$, with $k=1, \ldots, N_j$.
We observe that if the local fits $R_j$, $j=1,\ldots,d$, satisfy the interpolation conditions then the global PU approximant inherits the interpolation property. This trivially follows from the fact that the functions $W_j$, $j=1, \ldots, d$, form a partition of unity.

The coefficients $\{ c_k^j \}_{k=1}^{N_j}$ in \eqref{rad1} are determined by imposing the interpolation conditions
\begin{eqnarray}
\label{condinterp} 
R_j (\boldsymbol{x}_i^j )=f_i^j, \quad i=1,\ldots,N_j,
\end{eqnarray} 
which lead to the problem of solving $d$ linear systems of the form
\begin{equation}
A_j \boldsymbol{c}_j= \boldsymbol{f}_j ,
\label{sys1}
\end{equation}
where $  \boldsymbol{c}_j=  (c_1^j, \ldots, c_{N_j}^j  )^T$, $  \boldsymbol{f}_j = (f_1^j, \ldots , f_{N_j}^j  )^T$ and $A_j$ is the  local interpolation matrix whose entries are given by
\begin{equation}
(A_j)_{ik}= \phi ( || \boldsymbol{x}^j_i - \boldsymbol{x}_k^j ||_2 ), \quad i,k=1, \ldots, N_j.
\label{A}
\end{equation}
The system \eqref{sys1} admits a unique solution if and only if the matrix $A_j$ is non-singular. Micchelli gave conditions on $\phi$ which guarantee the non-singularity of $A_j$ \cite{Micchelli86}.  In particular, these conditions are fulfilled if $\phi$ is a strictly positive definite RBF. Nevertheless, even if here for simplicity we only consider strictly positive definite functions, we point out that such conditions are more general. Precisely, if $\phi$ is strictly conditionally positive definite of order $L$, a unique solution to the interpolation problem is ensured by adding to the expansion \eqref{rad1}  certain polynomials which form a basis for the $l$-dimensional  space $ \Pi_{L-1}^{M}$  of polynomials of total degree less than or  equal to $L-1$ in $M$ variables, where
	\begin{equation*}
	l=
	\begin{pmatrix} 
	L -1+M \\ 
	L-1 
	\end{pmatrix}.
	\end{equation*}

Among a large variety of known RBFs, which are defined in function of a shape parameter, we restrict our attention on strictly positive definite functions. Furthermore, we can distinguish between compactly  and globally supported RBFs. As examples of these two classes we consider the compactly supported Wendland $C^2$ and $C^6$ functions  and the  globally defined Mat\'ern $C^2$ and  Inverse MultiQuadric (IMQ) functions \cite{Fasshauer}. The latter is infinitely smooth and its analytic expression is
\begin{equation}
\phi(r)=\left(1+(\varepsilon r)^2 \right)^{-1/2} ,
\label{IMQ}
\end{equation}
where $r$ is the Euclidean norm and $\varepsilon$ is a positive shape parameter governing the flatness of the RBF. The family of Mat\'ern functions is instead characterized by a finite regularity. As example, the Mat\'ern $C^2$ is defined as
\begin{equation}
\phi (r) = e^{- \varepsilon r} (1+\varepsilon r).
\label{matern}
\end{equation}

Concerning Compactly Supported RBFs (CSRBFs), a well-known class of  functions which are smooth, locally supported and strictly positive definite is the family of Wendland's functions.
For instance, the  Wendland $C^2$ and $C^6$ functions are respectively defined as
\begin{equation}
\phi(r) = \left( 1-\varepsilon r\right)_+^{4} \left(4\varepsilon r+1\right),     
\label{WEN2}
\end{equation}
\begin{equation}
\phi(r) = \left( 1-\varepsilon  r\right)_+^{8} \left(32(\varepsilon  r)^3+25(\varepsilon  r)^2 + 8\varepsilon  r+1\right),   
\label{WEN}
\end{equation}
where  $( \cdot )_{+}$ denotes the truncated power function. 

Now, in order to give error bounds, we define two common indicators of data regularity: 
\begin{definition}
	The separation distance is given by
	\begin{equation}
	q_{ {\cal X}_N} = \frac{1}{2} \min_{  i \neq k} \left\| \boldsymbol{x}_i - \boldsymbol{x}_k \right\|_2.
	\label{sd}
	\end{equation}
\end{definition}
The quantity $ q_{ {\cal X}_N} $ represents the radius of the largest ball that can be centered at every point in ${\cal X}_N$ such that no two balls overlap.
\begin{definition}
	The fill distance, which is a measure of data distribution, is given by
	\begin{equation}
	h_{ {\cal X}_N, \Omega} =  \sup_{  \boldsymbol{x} \in \Omega} \left( \min_{ \boldsymbol{x}_k  \in {\cal X}_N} \left\| \boldsymbol{x} - \boldsymbol{x}_k \right\|_2 \right).
	\label{fd}
	\end{equation}
\end{definition}

\begin{remark}
	The distances \eqref{sd} and \eqref{fd}  give an idea of the node distribution, i.e. how uniform data are. Indeed, a set of data is supposed to be quasi-uniform with respect to a constant $C_{qu}$ if 
	\begin{equation}
	q_{ {\cal X}_N} \leq h_{ {\cal X}_N, \Omega} \leq C_{qu} q_{ {\cal X}_N}.
	\label{fill_d}
	\end{equation}
	More specifically, the definition of quasi-uniform points has to be seen in the context of more than one data set. The idea is to consider a sequence of such sets so that the domain $\Omega$ is more and more filled out. Then, points are said to be quasi-uniform if \eqref{fill_d} is satisfied by all the sets in this sequence with the same constant $C_{qu}$ \cite{Wendland05}.
	\label{fill}
\end{remark}

Furthermore, we need some technical considerations  on the regularity of the covering $ \{ \Omega_j \}_{j=1}^{d}$ and thus we give the following definitions.

\begin{definition}
	A subdomain $ \Omega_j \subseteq \mathbb{R}^M$ satisfies an interior cone condition if there exists an angle $ \theta \in (0,  \pi / 2)$ and a radius $\gamma>0$ such that, for all $ \boldsymbol{x} \in \Omega_j$, a unit vector $ \boldsymbol{\xi} ( \boldsymbol{x}) $ exists such that the cone
	\begin{equation*}
	C(\boldsymbol{x},\boldsymbol{\xi}(\boldsymbol{x}), \theta, \gamma)= \{ \boldsymbol{x} + \lambda \boldsymbol{y} : \boldsymbol{y}  \in  \mathbb{R}^M , || \boldsymbol{y} ||_2=1, \boldsymbol{y}^{T} \boldsymbol{\xi} ( \boldsymbol{x})  \geq cos( \theta) , \lambda \in  [0,\gamma]  \},
	\end{equation*}
	is contained in $ \Omega_j$.
\end{definition}

\begin{definition}
	Suppose that $ \Omega \subseteq  \mathbb{R}^M$ is bounded and $ {\cal X}_N= \{ \boldsymbol{x}_i  ,i=1, \ldots, N \} \subseteq \Omega$ is given. An open and bounded covering $ \{ \Omega_j \}_{j=1}^{d}$ is called regular for $( \Omega, {\cal X}_N)$, if the following properties are satisfied:
	\begin{itemize}
		\item[i.] for each $ \boldsymbol{ x} \in \Omega$, the number of subdomains $ \Omega_j$, with $ \boldsymbol{x} \in \Omega_j$ is  bounded by a global constant $C$,
		\item[ii.] each subdomain $ \Omega_j$ satisfies an interior cone condition,
		\item[iii.] the local fill distances $ h_{ {\cal X}_{N_j}, \Omega_j}$ are uniformly bounded by the global fill distance $h_{{\cal X}_N, \Omega}$.
	\end{itemize}
	\label{constant}
\end{definition}

Letting $C_{ \nu}^{k}  ( \mathbb{R}^{M} ) $ the space of all functions $f \in C^k$ whose derivatives of order $ | \boldsymbol{\beta} |=k $ satisfy $ D^{ \boldsymbol{\beta}} f ( \boldsymbol{x} ) = {\cal O} ( || \boldsymbol{x} ||_2^{ \nu} ) $ for $ || \boldsymbol{x} ||_2 \longrightarrow 0$, we consider the following convergence result \cite{Fasshauer,Wendland05}:

	\begin{theorem}
		Let $\Omega \subseteq  \mathbb{R}^M$ be open and bounded and assume that ${\cal X}_N = \{\boldsymbol{x}_i, i=1,$ $\ldots,N \}\subseteq \Omega$. Let $\phi \in C_{\nu}^k( \mathbb{R}^M)$ be a strictly conditionally positive definite function  of order $L$. Let $\{\Omega_j\}_{j=1}^{d}$ be a regular covering for $(\Omega, {\cal X}_N)$ and let $\{W_j\}_{j=1}^{d}$ be $k$-stable for $\{\Omega_j\}_{j=1}^{d}$. Then the error between $f \in {\cal N}_{\phi}(\Omega)$, where ${\cal N}_{\phi}$ is the native space of $\phi$, and its PU interpolant  is bounded by
		\begin{equation}
		| D^{\beta}f(\boldsymbol{x}) - D^{\beta}{\cal I}(\boldsymbol{x})  | \leq C^{'} h_{{\cal X}_N, \Omega}^{(k+\nu)/2 - |\beta|} |f|_{{\cal N}_{\phi}(\Omega)}, \nonumber
		\end{equation}
		for all $\boldsymbol{x} \in \Omega$ and all $|\beta| \leq k/2$. 
		\label{th1}
	\end{theorem}  
	
		\begin{remark}
			The first assumption in Definition \ref{constant} plays a crucial role also in the implementation of the PU method. In fact, such property leads to the requirement that the number of subdomains is proportional to the number of data  \cite{Wendland05}. 
			\label{remark_a}
		\end{remark}

	\subsection{Choosing suitable shape parameters and PU \\ subdomain sizes}
	\label{opt_ep}
	
	Usually,  the shape parameter $\varepsilon$ can greatly affect the accuracy of the resulting interpolant. Therefore, techniques allowing to select a predicted optimal shape parameter via error estimates have already been designed. Precisely, if the function is supposed to be known, the error can be exactly evaluated and thus the optimal shape parameter can be found without uncertainty. Otherwise, all the techniques based on error estimates give an approximated optimal value. Anyway, with abuse of notation, in what follows we will use the term optimal in the sense that  such approximation of the optimal value is \emph{close} to the one that can be found via trials and errors, for which the knowledge of the exact solution is supposed to be provided \cite{Fasshauer}.
	
	In the context of the PU method, aside from the value of the shape parameter, the size of the PU subdomains also plays a crucial role, especially when data with highly varying densities are considered. In literature, subdomains often consist of hyperspherical patches of the same radius $\delta$ \cite{Cavoretto15b,Fasshauer,Shcherbakov}.  Here, always considering  hyperspherical patches, we propose a novel method that allows to  suitably select  both the radius  $\delta_j$ and the shape parameter $\varepsilon_j$ for each PU subdomain $\Omega_j$, basing our considerations on an a priori error estimate. 
	
	We will focus on the so-called cross-validation algorithm, see \cite{Fasshauer,Golberg}, properly modified for a bivariate optimization problem. The cross-validation scheme has been firstly introduced in \cite{Allen,Golub79}. A variant of such method, known in literature as LOOCV, is detailed in \cite{Rippa}. Recent modifications of the cross-validation method can be found in \cite{FasshauerZhang}, where LOOCV is interpreted in the context of PDEs, and in \cite{TrahanWyatt}.
	
	Such approaches are always used in order to find the optimal value of the shape parameter for a global interpolation problem.  Here instead we are interested in selecting, for each PU subdomain, the optimal couple $(\delta_j,\varepsilon_j)$. 
	Carefully choosing, for each hypersherical patch such couple, leads to an accurate computation of the PU interpolant. In fact, supposing to have a regular covering, if we compare the result reported in Theorem \ref{th1} with the global error estimate shown in \cite{Wendland05}, we can see that the PU method preserves the local approximation order for the global fit. 
	Thus, the problem truly reduces in finding accurate local interpolants. In other words, if we improve the accuracy of the local fits, then also the one of the PU interpolant has benefits. This is even more evident from the following simple upper bound 
	\begin{equation*}
	\left|f(\boldsymbol{x})- {\cal I}(\boldsymbol{x}) \right| \leq \sum_{j=1}^{d} \left|f_j(\boldsymbol{x}) - R_j(\boldsymbol{x}) \right| W_j (\boldsymbol{x}) \leq \max_{j=1, \ldots, d} \left\|f_j-R_j \right\|_{L_{\infty}(\Omega_j)}, 
	\end{equation*}
	which shows that the PU approximation error is governed by the worst local error.
	
	Let us consider an interpolation problem on $\Omega_j$ of the form \eqref{rad1} and, for a fixed $i \in \{1, \ldots, N_j\}$, let
	\begin{equation*}
	R^{(i)}_j  ( \boldsymbol{x} )=  \sum_{k=1, k \neq i}^{N_j} c_k^j \phi ( || \boldsymbol{x} -\boldsymbol{x}^j_k ||_2 ),
	\end{equation*}
	be the $j$-th interpolant obtained leaving out the $i$-th data on $\Omega_j$. Moreover let
	\begin{equation}
	e^j_i= f^j_i-R^{(i)}_j (\boldsymbol{x}^j_i),
	\label{eq:er}
	\end{equation}
	be the  error at the $i$-th point. Then the quality of the local fit is determined by \emph{some} norm of the vector of errors $\boldsymbol{e}_j= ( e^j_1, \ldots ,e^j_{N_j})^{T},$ obtained by removing in turn one of the data points and comparing the resulting fit with the known value at the removed point.  Following \cite{Fasshauer,Rippa}, we can simplify the computation to a single formula by calculating
	\begin{equation}
	e^j_i= \dfrac{c^j_i}{\left(A_j \right)_{ii}^{-1}},
	\label{er0}
	\end{equation}
	where $c^j_i$ is the $i$-th coefficient of the RBF interpolant $R_j$ based on the full data set and $\left(A_j \right)_{ii}^{-1}$ is the $i$-th diagonal element of the inverse of the corresponding local interpolation matrix. 
	
	Precisely, in order to obtain an error estimate, we compute the following vector 
	\begin{equation}
	\left( e^j_1, \ldots, e^j_{N_j} \right)= \left( \dfrac{c^j_1}{ \left(A_j \right)_{11}^{-1}}, \ldots, \dfrac{c^j_{N_j }}{\left(A_j \right)_{N_j N_j}^{-1}} \right).
	\label{errori}
	\end{equation}
	In order to select the optimal couple $(\delta_j, \varepsilon_j)$ for each PU subdomain,  we compute \eqref{errori} for  several values of the radius  $(\delta_{j_1}, \ldots, \delta_{j_P})$ and of the shape  parameter $(\varepsilon_{j_1}, \ldots, \varepsilon_{j_Q})$.
	
	In \eqref{errori}, the dependence of the errors from the cardinality of the PU subdomain $N_j$, i.e. from the PU radius, is evident. Moreover, in this work, to stress the dependence of \eqref{errori} also from the shape parameter, for a fixed  $p \in \{1, \ldots, P \} $ and a fixed $q \in \{1, \ldots, Q \} $, we will use the notation 
	\begin{equation*}
	\boldsymbol{e}_j \left(\delta_{j_p},\varepsilon_{j_q} \right) = \left( e_1^j \left( \delta_{j_p},\varepsilon_{j_q} \right), \ldots, e^j_{N_j} \left( \delta_{j_p},\varepsilon_{j_q} \right) \right).
	\end{equation*}
	Thus, focusing on the maximum norm, we compute
	\begin{equation}
	E_j=
	\begin{pmatrix}
	||\boldsymbol{e}_j (\delta_{j_1},\varepsilon_{j_1}) ||_{\infty} & \cdots & ||\boldsymbol{e}_j (\delta_{j_1},\varepsilon_{j_Q}) ||_{\infty} \\	
	\vdots & \ddots & \vdots \\
	||\boldsymbol{e}_j (\delta_{j_P},\varepsilon_{j_1}) ||_{\infty} & \cdots & ||\boldsymbol{e}_j (\delta_{j_P},\varepsilon_{j_Q}) ||_{\infty} \\	
	\end{pmatrix}. 
	\label{mate}
	\end{equation}
	Note that \eqref{mate} provides an error estimate for several values of the PU radius and of the shape parameter. Therefore the $j$-th local approximant is computed  considering the couple $(\delta_j, \varepsilon_j)$ if
		\begin{equation}
		||\boldsymbol{e}_j(\delta_j, \varepsilon_j)||_{\infty}  = \min_{ p=1, \ldots, P} \left( \min_{  q=1, \ldots, Q} (E_j)_{pq} \right).
		\label{opt_coppia}
		\end{equation}
	In other words, the BLOOCV-PU interpolant assumes the form
	\begin{equation}
	{\cal \tilde{I}}( \boldsymbol{x})= \sum_{j=1}^{d} \tilde{R}_j( \boldsymbol{x} ) W_j ( \boldsymbol{x}), \quad \boldsymbol{x} \in \Omega,
	\label{intg_BLOOCV}
	\end{equation}
	where, for each subdomain $\Omega_j$, $\tilde{R}_j$ is given by
	\begin{eqnarray}
	\label{rad1_BLOOCV}
	\tilde{R}_j({\boldsymbol{x}})=\sum_{k=1}^{\tilde{N}_j} c^j_{k} \phi_{\varepsilon_j}  (|| \boldsymbol{x}-\boldsymbol{x}^j_k ||_2 ),
	\end{eqnarray}
	and $\tilde{N}_j$ indicates the number of points in $\Omega_j$ of radius $\delta_j$. 
	
	Observe that, consistently with Definition \ref{constant} and Remark \ref{remark_a},  if the number of patches is proportional to $N$ and if the subdomains form a covering of $\Omega$, then such covering is also regular. This trivially follows from the fact that a hypersphere of radius $\delta_j$ always satisfies an interior cone condition with constants independent from the space dimension; precisely, $\gamma= \delta_j$ and $\theta= \pi /3$ \cite{Wendland05}. Therefore, all the considerations made in the previous subsection also hold for the BLOOCV-PU interpolant. 
	
	 This approach obviously leads to a benefit in terms of accuracy, especially when irregular data are considered. However, we have to point out that the computation of the inverse for each couple  $(\delta_{j_p},\varepsilon_{j_q})$ is particularly costly  for large $\delta_{j_p}$. Therefore, for each $\Omega_j$ we need to carefully choose the extreme values of the discrete searching range for the radius, i.e. the interval $[\delta_{j_1}, \delta_{j_P}]$.
	
	 Precisely, for each PU subdomain,  we have at first to fix the  intervals $[\delta_{j_1}, \delta_{j_P}]$ and $[\varepsilon_{j_1},  \varepsilon_{j_Q}]$, used to find out  $(\delta_j,\varepsilon_j)$. Many researchers already worked on the problem of finding suitable values for shape parameter in order to increase the accuracy and, at the same time, avoid problems of instability. Thus, one can easily guess how to select a good range for the shape parameter \cite{Demarchi15,Driscoll-Fornberg02,Fasshauer,Fasshauer15,Fornberg11}. In other words, for what concerns the shape parameter the notation simplifies, since for each subdomain we can consider the same  discrete values, namely  $(\varepsilon_{1}, \ldots, \varepsilon_{Q})$. On the opposite, fixing for all the subdomains the same  discretization $(\delta_{1}, \ldots, \delta_{P})$ can lead to inaccurate solutions. More specifically:
	\begin{problem} Arbitrarily fixing, for all the subdomains,  the same searching interval   $[\delta_{1},  \delta_{P}]$ can lead to the following  issues:
		\begin{itemize}
			\item[1.] the union of the PU subdomains might not form a covering of the domain;
			\item[2.] in regions characterized by a  low density of points the interval  $[\delta_{1},  \delta_{P}]$  can be too small to avoid empty patches or subdomains containing very few points;
			\item[3.] in regions characterized by a  high density of points, the interval  $[\delta_{1},  \delta_{P}]$  can be too large and, in this case, both complexity and ill-conditioning grow.
		\end{itemize}
		\label{problem1}
	\end{problem}
	
	In the subsequent section we detail a feasible scheme useful to determine the  interval  $[\delta_{1},  \delta_{P}]$, in which we can search for the $j$-th suitable radius, avoiding the above mentioned problems. To reach this aim, we first need an efficient partitioning structure, used to organize points among the different subdomains. Therefore, we propose a novel multidimensional procedure, built ad hoc for the PU method and independent from the problem geometry.
	
	\section{Feasible computation of the BLOOCV-PU \\interpolant}
	\label{comput_asp}
	
	As already pointed out in Subsection \ref{opt_ep}, the searching interval  $[\delta_{1},  \delta_{P}]$  must be properly selected. Essentially, in order to obtain reliable error estimates, we want to make sure of having a reasonable number of points on each patch. Such consideration suggests the use of a $K$-nearest neighbor procedure. As example in \cite{Deparis}, suitable supports of CSRBFs  have been selected  detecting, via a triangulation, the $K$-nearest neighbor set \cite{Arya98}. This turns out to be expensive and moreover fixing an arbitrary $K$ does not guarantee a good approximation. Thus, we will use a similar procedure to \cite{Deparis} only to determine the initial reasonable searching range  $[\delta_{1},  \delta_{P}]$  for the radius of the $j$-th patch and then such interval will be used in the computation of \eqref{mate}.
	
	The complexity needed to construct the BLOOCV-PU interpolant will be taken into account. Specifically, we will not perform a $K$-nearest neighbor procedure, but we will use the new multidimensional I-PS. It leads to a considerable saving in terms of computational time with respect to \cite{Cavoretto15b,AIP_15Emma} and unlike them it can be applied in any space dimension $M$ (and not only for $M=2,3$).
	
	We will treat the problem in the most general setting. Thus let us consider a set of data ${\cal X}_N=\{\boldsymbol{x}_i \in \Omega , i=1,\ldots,N\}$, where $\Omega \subseteq \mathbb{R}^M$ is a simply connected region. In order to perform the BLOOCV-PU method, we need to define an hyperrectangle ${\cal R}_M$ containing the scattered data
	\begin{align*}
	{\cal R}_{M}= \prod_{m=1}^{M}   \left[ \min_{i=1, \ldots, N}  x_{im}, \max_{i=1, \ldots, N}  x_{im} \right],
	\end{align*}
	and the bounding box containing the nodes, i.e. the box of edge
	\begin{equation}
	l_{box}= \max_{m=1,\ldots,M} \left( \max_{i=1, \ldots, N} x_{im} \right) - \min_{m=1,\ldots,M} \left( \min_{i=1, \ldots, N} x_{im} \right).
	\label{lato_box1}
	\end{equation}
	
	Then, consistently with Remark \ref{remark_a}, we define the PU centres as a grid of $ d^{{\cal R}_M} = \left( d_{PU}^{{\cal R}_M} \right)^M$ points on ${\cal R}_{M}$,  where
	\begin{align*}
	d^{{\cal R}_{M}}_{PU}=\bigg \lfloor  \frac{\displaystyle 1}{\displaystyle 2} l_{box} \left(\frac{N}{V_{\Omega}}\right)^{1/M}\bigg \rfloor, 
	\end{align*}
	and $V_{\Omega}$ is the hypervolume of the simply connected region $\Omega$. Then, to make sure that patches form a covering of the domain, we can set the radii of the hyperspheres $\delta_j$ such that
	\begin{equation}
	\delta_j \geq \frac{\displaystyle l_{box}}{\displaystyle d_{PU}^{{\cal R}_M}}, \quad j=1, \ldots, d.
	\label{raggio_pu}
	\end{equation}
	With \eqref{raggio_pu} we solve the  issue $1.$ outlined in Problem \ref{problem1}. 
	
	The initial number of subdomains $d^{{\cal R}_M} $ is later reduced by taking only those $d$ centres lying in $\Omega$. This step provides the set of PU centres ${\cal C}_d= \{ \boldsymbol{\bar{x}}_i, i=1, \ldots, d \}$, used to construct the PU interpolant. Note that, because of \eqref{raggio_pu}, the set ${\cal C}_d$ forms a covering for $\Omega$. Furthermore, in the same way we also define a  set ${\cal E}_s= \{ \boldsymbol{\tilde{x}}_i, i=1, \ldots, s \}$  on $\Omega$, which is used to evaluate the unknown function via BLOOCV-PU interpolation.
	Then, in order to organize points into the different patches and consequently choose a suitable searching interval for the radius (see Subsection \ref{select}), we consider the partitioning structure described in Subsection \ref{ips}.

	\subsection{Multidimensional integer based partitioning structure}
	\label{ips}
	To make simpler the presentation, we first consider hyperspherical patches all having the same radius \cite{AIP_15Emma}
	\begin{equation}
	\delta = \frac{\displaystyle l_{box}}{\displaystyle d_{PU}^{{\cal R}_M}}.
	\label{fisso}
	\end{equation}
	Then, in order  to solve the local interpolation problems, we need to develop a procedure enabling us to store the points among the different PU subdomains. Such scheme must be independent from the problem geometry and work in any dimension, as kd-trees \cite{Arya98,deBerg97,Fasshauer,Yaoa}. These are effective and widely used numerical tools, but they are not specifically constructed for the PU method. 
	
	Thus, starting from the bivariate and trivariate procedures, that in the following we will call the Sorting-based Partitioning Structures (S-PSs) \cite{Cavoretto15b,AIP_15Emma}, our aim is to build a multidimensional procedure which allows to consider variable radii $\delta_j$. As in the S-PS, we store the points  into $q^M$ blocks, where
	\begin{equation} 
	q= \bigg \lceil \frac{\displaystyle l_{box}}{\displaystyle \delta} \bigg \rceil.
	\label{q1}
	\end{equation}
	More precisely, we number blocks from $1$ to $q^M$, starting from the subspace of dimension $M-1$, obtained projecting along the first coordinate and thus parallel to the remaining ones. In order to fix the idea, in a 2D context they are numbered from bottom to top, left to right. 
	
	\begin{remark}
		In bivariate interpolation blocks are generated by the intersection of two orthogonal strips. In multivariate problems blocks are generated by the intersection of $M$ hyperrectangles. In what follows with abuse of notation we will continue to call such hyperrectangles with the term  strips.
	\end{remark}
	
	Then, in order to store the points among the different patches the following computational issue, known as \emph{containing query}, needs to be solved
	\begin{itemize}
		\item given a PU centre $\boldsymbol{\bar{x}}_j$, find the $k$-th block containing the centre.
	\end{itemize}
	Such problem can be easily solved taking into account that, given a PU centre $\boldsymbol{\bar{x}}_j$, if $k_m$ is the index of the strip parallel to the subspace of dimension $M-1$ generated by $x_r$, $r=1, \ldots,M$ and $r \neq m$, containing the $m$-th coordinate of $\boldsymbol{\bar{x}}_j$, then the index  of the  $k$-th block  containing the subdomain centre is
	\begin{align} 
	k=\sum_{m=1}^{M-1} \left( k_m-1 \right) q^{M-m}+k_M.
	\label{idx_2}
	\end{align}
	To find the indices $k_m$, $m=1, \ldots, M$, in \eqref{idx_2}, we use an integer-based procedure consisting in rounding off to an integer value. Specifically, for each PU centre $\boldsymbol{\bar{x}}_j=(\bar{x}_{j1}, \ldots ,\bar{x}_{jM})$, we have that
	\begin{equation}
	k_m= \bigg\lceil \frac{\bar{x}_{jm}}{\delta} \bigg\rceil.
	\label{idx_3}
	\end{equation}
	
	Then, exactly the same procedure is adopted in order to store into the different blocks both scattered data and evaluation points, i.e.  the I-PS assigns:
	\begin{itemize}
		\item[i)] to each scattered point $\boldsymbol{x}_i$, $i=1, \ldots ,N$, the index of the block in which it lies,
		\item[ii)] to each evaluation point $\boldsymbol{\tilde{x}}_i$, $i=1, \ldots ,s$, the index of the block in which it lies.
	\end{itemize}
	
	Moreover, always supposing to have a fixed radius and assuming that the $j$-th centre belongs to the $k$-th block, from \eqref{q1} the fact that we search for the points lying in the $j$-th patch among those lying in the $k$-th block and in its $3^M-1$ neighboring blocks easily follows. 
	
	On the opposite, here the radius is supposed to be variable for each patch and thus, if the radius $\delta_j$ is such that
	\begin{equation*}
	\delta_j > n \delta, \quad n \in \mathbb{N}^{+},
	\end{equation*}
	given the centre $\boldsymbol{\bar{x}}_j$ we search for the neighboring points in the $k$-th block and in its $(3+2^n)^M-1$ neighboring blocks. 
	
	In \cite{Cavoretto15b,AIP_15Emma} nodes and evaluation points are organized in blocks by using recursive calls to a sorting routine, while here this step is replaced by \eqref{idx_2} and \eqref{idx_3}. Such approach enables us to improve three aspects of the partitioning structure presented in the last mentioned papers. Precisely, the I-PS:
	\begin{itemize}
		\item[1)] works in any dimension,  while the S-PS only works for $M=2,3$,
		\item[2)] allows to work with variable radii, while the  S-PS strictly depends on a fixed size of the subdomains,
		\item[3)] reduces the complexity of the sorting-based storing procedure (see Section \ref{analisi_costi}). 
	\end{itemize}
	
	\subsection{Selection of a  searching interval for the PU radius}
	\label{select}
	
	When we deal with quasi-uniform  or grid data, the number of points in each subdomain of radius $\delta$ is about constant. On the opposite, in case of  irregular nodes,  the number of points lying in the different patches is far from being constant or, even worst, we can have empty subdomains. This consideration turns out to be useful to determine the lower bound of the interval  $[\delta_{1},  \delta_{P}]$ for the $j$-th subdomain. In fact, given  the number $N$ of scattered data  in $\Omega$  and  its hypervolume $V_{\Omega}$,  from a simple proportion we have that a suitable number $K$ of points belonging to $\Omega_j$ of radius \eqref{raggio_pu} is
	\begin{equation}
	K \approx \frac{N B(\delta)}{V_{\Omega}},
	\label{punti_vic}
	\end{equation}
	where $B(\delta)$ is the hypervolume of the hypersphere of radius $\delta$, defined as in \eqref{fisso}. 
	The value found in \eqref{punti_vic} represents the number of points we expect on average on each  patch supposing to have a uniform node distribution.
	
	Therefore, given the $j$-th subdomain of radius $\delta_{j_1}=\delta$, we compute its cardinality via the I-PS. Then, if such cardinality  is less than the one given by \eqref{punti_vic}, $\delta_{j_1}$ is updated as follows
	\begin{equation}
	\delta_{j_1} = \delta_{j_1}+t \delta,
	\label{raggio_pu1}
	\end{equation}
	where $0<t<1$. The procedure continues in this way until \eqref{punti_vic} is satisfied, i.e. $\delta_{j_1}$ is determined, with recursive calls to the I-SP, so that
	\begin{equation}
	Card(\Omega_j)  \geq \frac{N B(\delta_{j_1})}{V_{\Omega}}.
	\label{card_omega}
	\end{equation}
	Acting in this way, we solve the computational issue  $2.$ outlined in Problem \ref{problem1},  i.e. there are few enough points for each patch.
	
	Then, in order to avoid also the third issue of Problem \ref{problem1}, the simplest strategy, which takes into account the density of points and turns out to be effective, is to choose  $P$ discrete values in an interval of the form 
	\begin{equation}
	 [\delta_{j_1}, h \delta_{j_1}], \quad h \in \mathbb{R}^{+},  \quad h>1.
	\label{range}
	\end{equation}
	Roughly speaking, since the upper bound of the searching interval is proportional to the lower bound and since this lower bound is large only if the density of points is low, we effectively avoid problems arising from high density of points, i.e. systems are not too large and the ill-conditioning is kept under control.
	
\begin{figure}[ht!] 
				\begin{center}
					\makebox[\textwidth]{
			\includegraphics[height=.28\textheight]{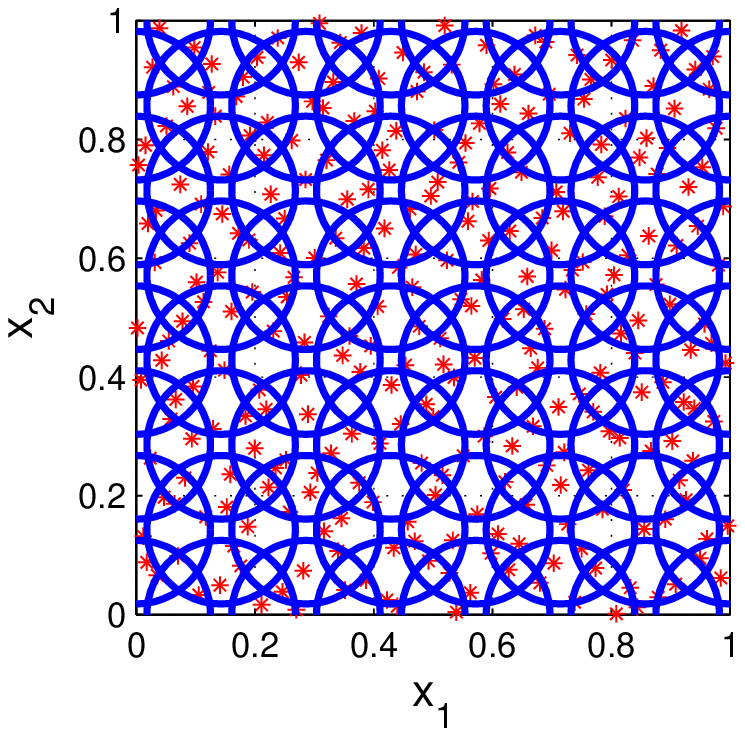} 
			\hskip -1.5cm
			\includegraphics[height=.28\textheight]{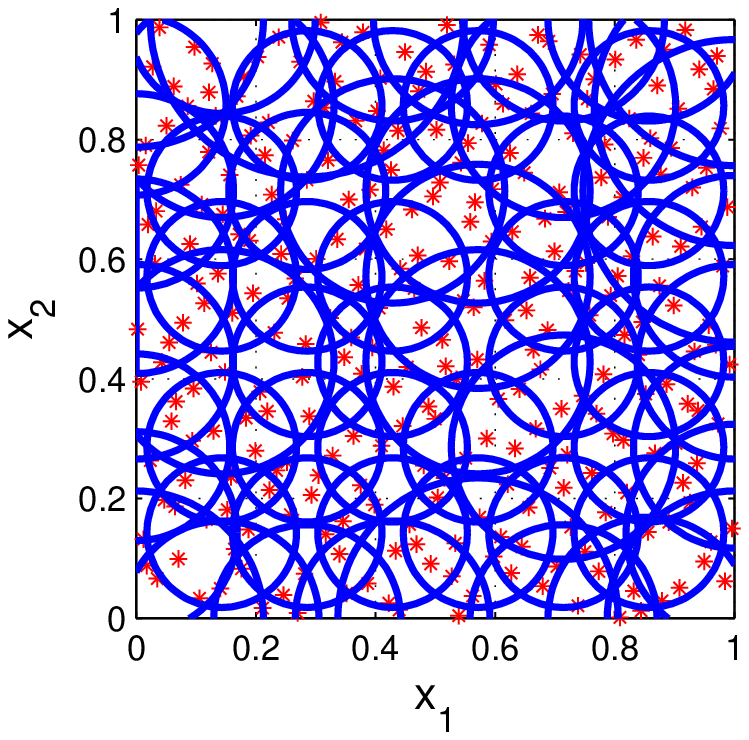} }
			\caption{Examples of PU structures  covering $289$ Halton data in $\Omega=[0,1]^2$: the classical PU structure (left) and the variable structure obtained via the BLOOCV-PU algorithm (right).}
			\label{figuraPUM_var}
		\end{center}
	\end{figure}

	We end this section with  the illustrative Figure \ref{figuraPUM_var}, devoted to show how the classical PU structure is modified by means of the BLOOCV-PU algorithm. In the left frame we plot the classical structure by choosing $\delta_j=\delta$, $j=1, \ldots ,d$, while in the right frame we show the result of the BLOOCV-PU method.
	
	\section{Complexity analysis}
	\label{analisi_costi}
	In this section we analyze the computational cost of the proposed method. It will be pointed out that the main cost is due to the computation of the error estimates, while the I-PS turns out to be really cheap.
	
	\subsection{The BLOOCV local implementation}
	For each PU subdomain, several error estimates are calculated via \eqref{er0}. Such computation needs ${\cal O} (N_j^3)$ operations. Thus, simplifying the calculation via \eqref{er0} is the key step which enables us to maintain a reasonable complexity cost. Indeed, evaluating the error via \eqref{eq:er} is  computationally expensive. In particular, the matrix inverse, which requires ${\cal O} (N_j^3)$ operations, must be computed  for each node. This step needs a total computational cost of ${\cal O} (N_j^4)$ operations, $j=1, \ldots, d$, but using \eqref{er0}, the complexity cost significantly decreases. 
	
	 However, the complexity of the proposed algorithm is quite high. The error estimate \eqref{er0} needs to be computed for each subdomain $\Omega_j$ and for each $\delta_{j_i}$, $i=1, \ldots, P$, and $\varepsilon_{k}$, $k=1, \ldots, Q$.
	
	\begin{remark}
		The computation of the error estimate for the shape parameter  can be slightly speeded up  by using the \textsc{Matlab} routine
		{\tt CostEpsilon.m}, proposed in \cite{Fasshauer}. However,  the same approach cannot be performed for an optimal choice of the radius $\delta_j$. Therefore, in our free software package, for easiness of the \textsc{Matlab} user, we carried out the standard implementation as in \eqref{mate}. 
	\end{remark}

	\subsection{The integer-based partitioning structure}
	The I-PS, after organizing the scattered data into the different blocks, given a subdomain $\Omega_j$ searches for all the points lying in $\Omega_j$ in a reduced number of blocks. Specifically, in order to store the scattered data among the different blocks, it makes use of an integer-based procedure that assigns to each node $N_i$, $i=1,\ldots,N$, the corresponding block. This step requires ${\cal O} (N)$ time. Then, we apply the optimized searching routine already used in \cite{Cavoretto15b,AIP_15Emma}. Such procedure is performed in a constant time (refer to \cite{AIP_15Emma} for further details).
	
\begin{figure}
	\begin{center}
		\makebox[\textwidth]{
		\includegraphics[height=.28\textheight]{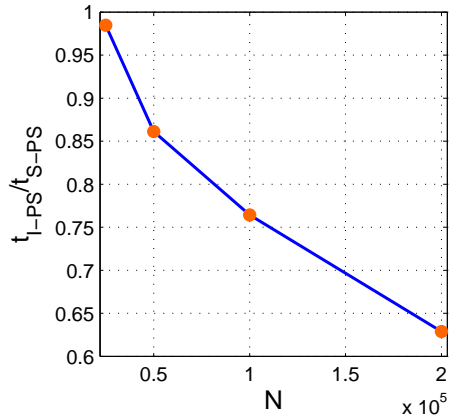}}
		\caption{CPU time ratios $t_{I-PS}/t_{S-PS}$ by varying $N$.} 
		\label{fig_tempi_ratio}
		\end{center}
	\end{figure}

	Observe that the I-PS turns out to be more efficient than the S-PS; in fact the latter, to store the points among the different blocks, needs ${\cal O} (N \log N)$ operations.
	Table \ref{tabe_tib} and Figure \ref{fig_tempi_ratio} support our findings. Specifically,  we consider in a 2D framework different sets of Halton data. Tests have been carried out on a Intel(R) Core(TM) i7 CPU 4712MQ 2.13 GHz processor.
	
	\begin{table}[h]
		\caption{CPU times (in seconds) obtained by running the sorting-based procedure ($t_{S-PS}$) and the integer-based  one ($t_{I-PS}$).}
		\begin{center}
			\begin{tabular}{ccccc} 
				\hline\noalign{\smallskip}
				$N$	& $25000$ & $50000$  & $100000$  & $200000$ \\
				\noalign{\smallskip}
				\hline
				\noalign{\smallskip}
				$t_{I-PS}$  & $5.13$ & $10.68$ & $21.99$ & $45.00$\\
				$t_{S-PS}$ & $5.21$ & $12.40$ &  $28.77$ & $71.55$\\
				\hline 
			\end{tabular}
		\end{center}
		\label{tabe_tib}
	\end{table}
	
	\subsection{Computation of the PU interpolant}
	The computation of the local interpolation  problems consists in solving $d$ linear systems of size $N_j \times N_j$, with $N_j \ll N$. This step involves   a computational cost of ${\cal O}(N^3_j)$ time, $j =1, \ldots , d$.
	Since the number $d$ of subdomains is bounded by ${\cal O}(N)$, this leads to ${\cal O}(N)$ operations for solving all of them. 
	Finally, in order to get the value of the global fit \eqref{intg}, we have to add up a constant number of local RBF interpolant. This requires ${\cal O}(1)$ time.

	\section{Numerical experiments}
	\label{Numerica_ex}
	
	This section is devoted to show, by means of extensive numerical simulations the flexibility and the accuracy of the proposed BLOOCV-PU method. It is applied fixing the initial intervals for the radii as in \eqref{range}, with $h=2$ and $P=6$. These values are chosen so that we ensure to have enough points on each subdomain and a sufficient number of radii to test the accuracy of the interpolants. Moreover, we will use as PU weights the Wendland $C^2$ function, see \eqref{WEN2}.  
	
	Tests are carried  out considering  the so-called \emph{product} and \emph{valley functions} \cite{Bozzini,Nielson}, respectively defined as:
	\begin{equation*}
	f_1(x_1,x_2)  = 16 x_1 x_2 (1-x_1) (1-x_2), \quad f_2(x_1,x_2) = \frac{1}{2}  x_2 \left[ \cos(4x_1^2+x_2^2-1) \right]^4 .
	\end{equation*}
	
	To  point out the accuracy of the BLOOCV-PU interpolant, we will refer to the Maximum Absolute Error (MAE) and  the Root Mean Square Error (RMSE), whose formulas are:
	\begin{eqnarray} \label{MAE}
	\textrm{MAE} = \max_{ i=1, \ldots, s} |f(\tilde{\boldsymbol{x}}_i) - {\cal \tilde{I}}(\tilde{\boldsymbol{x}}_i) |, \quad \textrm{RMSE} = \sqrt{\frac{1}{s}\sum_{i=1}^{s} |f(\tilde{\boldsymbol{x}}_i) - {\cal \tilde{I}}(\tilde{\boldsymbol{x}}_i)|^2},
	\end{eqnarray}
	where $\tilde{\boldsymbol{x}}_i$, $i=1, \ldots, s$, forms a grid of $40 \times 40$  points in which the interpolant is sought.
	
\begin{figure}
		\begin{center}
		\makebox[\textwidth]{
		\includegraphics[height=.28\textheight]{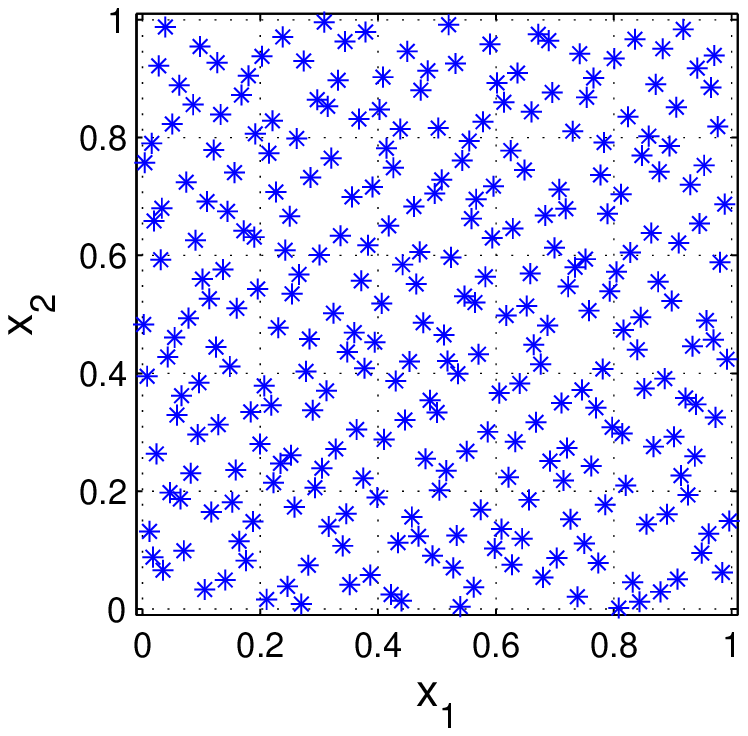}
		\hskip -1.5cm
		\includegraphics[height=.28\textheight]{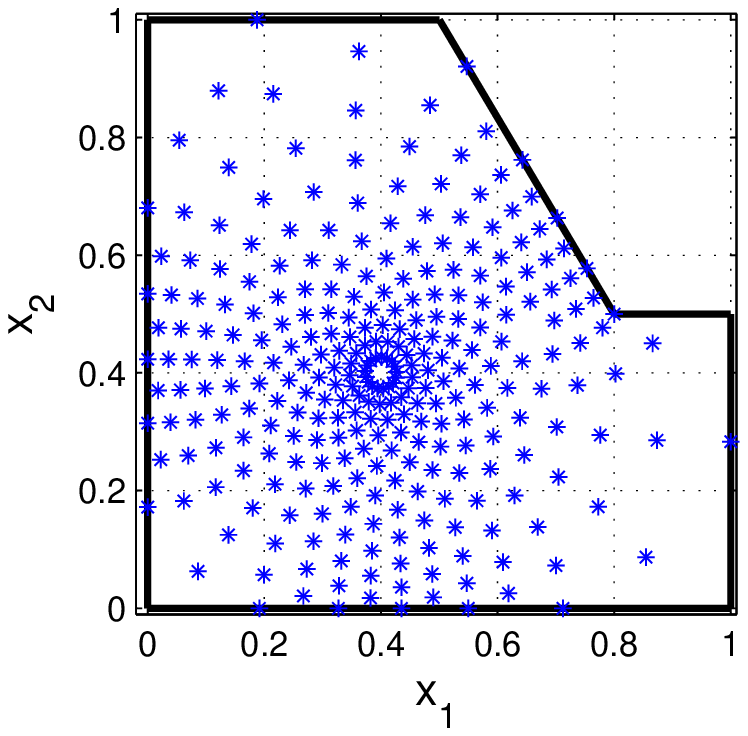}}\\
		\caption{Examples of data sets: $289$ Halton points (left) and  $289$ points in a polygonal region coming from a Schwarz-Christoffel transformation (right).}
		\label{concentric_circles_400_cerchio}
	\end{center}
	\end{figure}

	Concerning the data sets used in our numerical experiments we take, as quasi-uniform points, the Halton data \cite{Fasshauer}, see Figure \ref{concentric_circles_400_cerchio} (left). 
	Then, in order to test the method with particularly hard nodes, we consider points coming from a Schwarz-Christoffel transformation, refer to \cite{Driscoll,Driscoll02,Heryudono} for further details. More precisely, we focus on a special case of conformal map from the unit disk onto a polygon. Thus, at first we define nodes in the unit disk and then we map them into a chosen polygon. As example, in Figure \ref{concentric_circles_400_cerchio} (right) we show the result of  conformally mapping onto a simply connected region points on concentric circles. We consider such points because  they are far from being quasi-uniform,  same time they are constructed with a specific rule and thus tests are repeatable. Moreover, they simulate practical situations, as it will be evident in Section \ref{foresta} when we will deal with data coming from real life.
	
	The BLOOCV-PU approach will be compared with the classical PU method, i.e. the PU scheme is applied with a fixed radius and a fixed shape parameter \cite{Cavoretto15b,Fasshauer,Wendland05}.  Thus, such classical approach requires to fix these two parameters. Concerning the radius, the classical PU method is generally applied considering a fixed size of the patches as in \eqref{fisso}, while for the shape parameter, in literature, the choice is \emph{almost arbitrary}. In fact, even if techniques allowing to obtain stable approximations when $\varepsilon \rightarrow 0$ have already been developed \cite{Demarchi15,Fornberg11}, there is not an a priori good value for the latter. We will later discuss this concept providing several tests.
	
	In our numerical experiments, we expect the following behavior classes depending  on the distribution of the data set:
	\begin{itemize}
		\item[i.] with quasi-uniform points: BLOOCV-PU and classical PU both give accurate approximations,
		\item [ii.] with non-conformal points: the classical PU method fails, while the BLOOCV-PU maintains a good accuracy.
	\end{itemize}

	\subsection{Tests with quasi-uniform data}
\begin{figure}
		\begin{center}
			\makebox[\textwidth]{
		\includegraphics[height=.24\textheight]{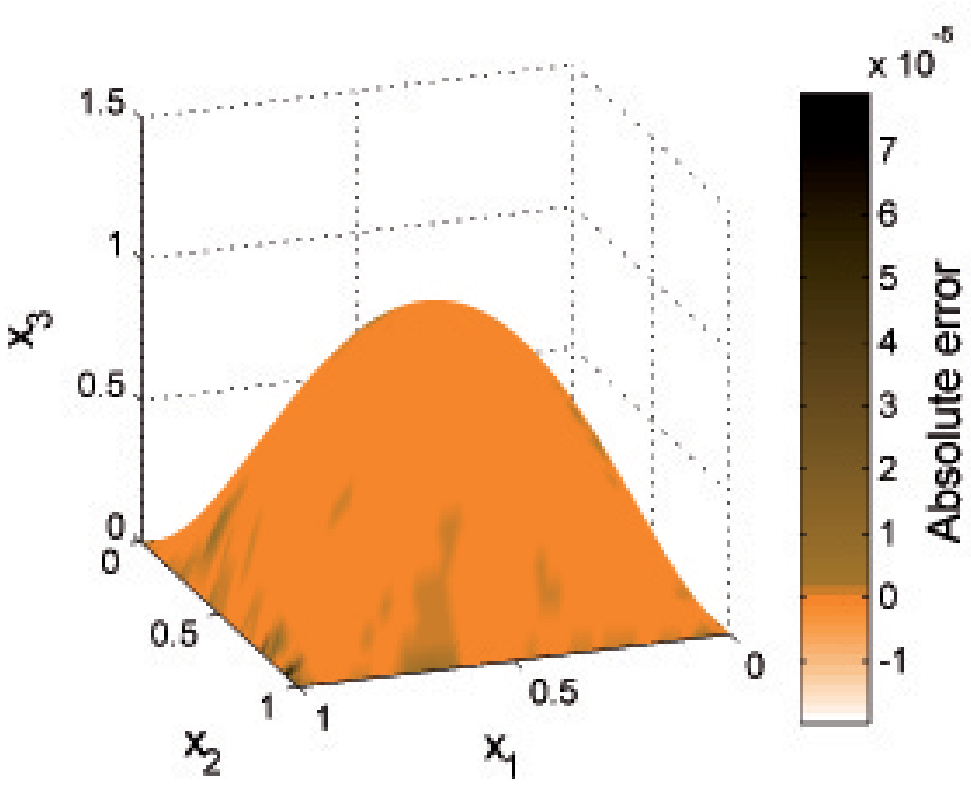}
		\hskip -0.2cm
		\includegraphics[height=.24\textheight]{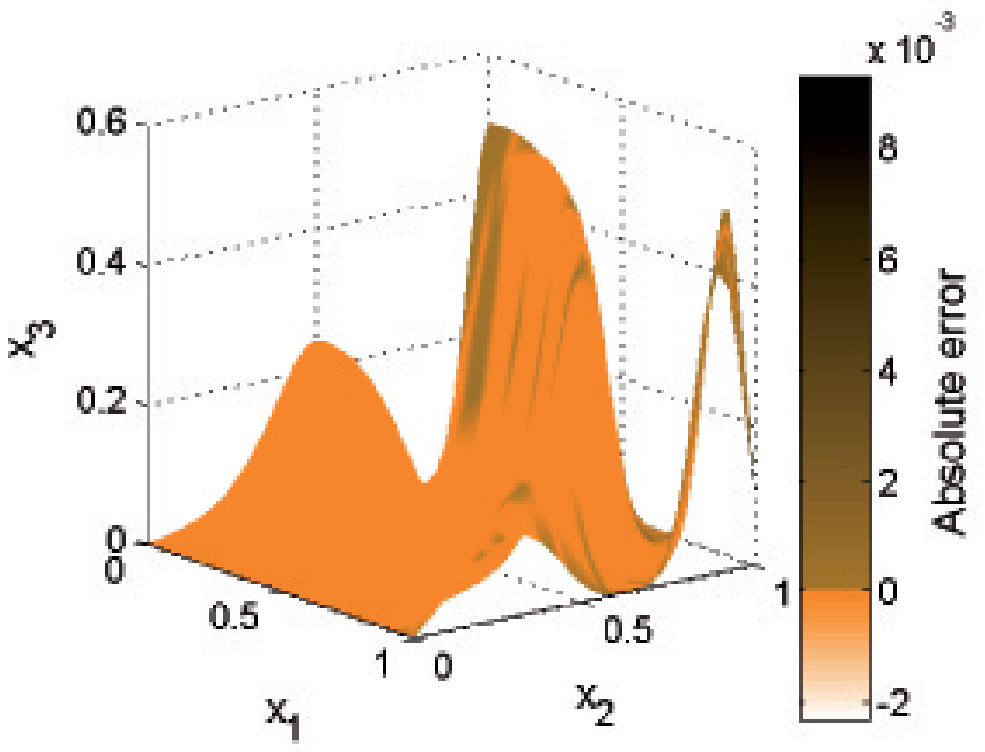}}\\
		\caption{The false-colored absolute errors computed on $1089$ Halton points and obtained by applying the BLOOCV-PU method with the IMQ as local RBF interpolant for $f_1$ (left) and $f_2$ (right).} 
		\label{fig_fun}
	\end{center}
	\end{figure}

	In this subsection we consider Halton points and, because of their regularity, we take a smooth RBF as local approximant, specifically the IMQ, see \eqref{IMQ}. Therefore, choosing $30$ values for the shape parameter in the  interval $[0.1, 10]$ is reasonable. Results, obtained by means of the BLOOCV-PU interpolant with the functions $f_1$ and $f_2$, are shown in Tables \ref{tabeu1} and \ref{tabeu2}, respectively. We also report  the errors of the classical PU method, obtained by fixing the radius as in \eqref{fisso} and the shape parameter $\varepsilon=0.6$. For a graphical representation of the distribution of the absolute error refer to Figure \ref{fig_fun}. Such figure shows that the error for  both test functions is larger close to the boundary. Moreover,  for $f_2$ it also increases in correspondence of the oscillations of the test function.
	
	 From Tables \ref{tabeu1} and \ref{tabeu2}, we can easily note that the results obtained with the BLOOCV-PU interpolant are more accurate than the ones carried out with the standard PU. Anyway, in order to get such accurate approximation, we have to pay in terms of efficiency. For instance, with $289$ and $1089$ points the BLOOCV-PU can be computed in $5.48$ s and $20.08$ s, respectively, while the classical PU interpolation only requires $0.21$ s and $0.60$ s.

	
	\begin{table}[ht!]
		\begin{center}
			\begin{tabular}{cccc} 	\hline\noalign{\smallskip}
				$N$ & method	 & RMSE & MAE  \\
				\hline 
				\rule[0mm]{0mm}{3ex}
				$\hskip-2pt 289$ & PU &  $3.64{\rm E}-03$    &  $5.66{\rm E}-02$ 	   \\
				&  BLOOCV-PU          &  $1.03{\rm E}-05$    &  $2.36{\rm E}-04$  	   \\
				\rule[0mm]{0mm}{3ex}
				$1089$     &  PU     & $7.57{\rm E}-04$    &  $1.52{\rm E}-02$	   \\
				&   BLOOCV-PU  & $2.88{\rm E}-06$    &  $7.89{\rm E}-05$	  \\
				\rule[0mm]{0mm}{3ex}
				$4225$     & PU   & $3.88{\rm E}-04$    &  $1.01{\rm E}-02$ 	   \\
				&   BLOOCV-PU   & $3.84{\rm E}-07$    &  $1.39{\rm E}-05$ 	      \\
				\rule[0mm]{0mm}{3ex}
				$16641$     &  PU   & $8.27{\rm E}-04$    &  $3.27{\rm E}-02$   \\
				&   BLOOCV-PU   & $9.67{\rm E}-08$    &  $3.15{\rm E}-06$      \\
				\rule[0mm]{0mm}{3ex}
				$66049$     &  PU  & $1.08{\rm E}-05$    &  $1.09{\rm E}-04$   \\
				&  BLOOCV-PU     & $2.68{\rm E}-08$    &  $6.80{\rm E}-07$      \\
				\hline 
			\end{tabular}
		\end{center}
		\caption{RMSEs and MAEs computed on Halton points and obtained by using the IMQ as local RBF interpolant for $f_1$.}
		\label{tabeu1}
	\end{table}

	\begin{table}[ht!]
		\begin{center}
			\begin{tabular}{cccc} 	\hline\noalign{\smallskip}
				$N$ & method	 & RMSE & MAE  \\
				\hline 
				\rule[0mm]{0mm}{3ex}
				$\hskip-2pt 289$ & PU  & $2.59{\rm E}-02$    &  $4.30{\rm E}-01$	   \\	
				&  BLOOCV-PU         & $1.32{\rm E}-02$    &  $2.76{\rm E}-01$   	   \\
				\rule[0mm]{0mm}{3ex}
				$1089$     &  PU    & $3.51{\rm E}-03$    &  $6.20{\rm E}-02$   \\
				&   BLOOCV-PU  & $2.11{\rm E}-04$    &  $8.93{\rm E}-03$   \\
				\rule[0mm]{0mm}{3ex}
				$4225$     & PU   & $8.63{\rm E}-04$    &  $2.00{\rm E}-02$	   \\
				&   BLOOCV-PU    & $3.88{\rm E}-06$    &  $1.12{\rm E}-04$ 	      \\
				\rule[0mm]{0mm}{3ex}
				$16641$     &  PU   & $4.07{\rm E}-04$    &  $1.18{\rm E}-02$   \\
				&   BLOOCV-PU   & $8.26{\rm E}-08$    &  $2.80{\rm E}-06$        \\
				\rule[0mm]{0mm}{3ex}
				$66049$     &  PU   & $1.23{\rm E}-04$    &  $4.19{\rm E}-03$     \\
				&  BLOOCV-PU    & $5.10{\rm E}-08$    &  $1.76{\rm E}-06$     \\
				\hline 
			\end{tabular}
		\end{center}
		\caption{RMSEs and MAEs computed on Halton points and obtained by using the IMQ as local RBF interpolant for $f_2$.}
		\label{tabeu2}
	\end{table}

	One may argue that the results of the classical PU method might truly improve by varying the shape parameter. This is trivially true, but at the same time  automatically choosing a \emph{safe} shape parameter, which gives reliable approximations, is one of the main advantages of the BLOOCV-PU method. Anyway, in order to clarify this concept, we report in Table \ref{tabe_1a} the results obtained by choosing in the classical PU algorithm the optimal shape parameter, but not the optimal radius, for each subdomain. Note that this is only a particular case of the BLOOCV-PU algorithm and thus no longer investigated.
	
	\begin{table}  
		\begin{center}
			\begin{tabular}{ccc}
				\hline\noalign{\smallskip}
				N  &  RMSE &  MAE  \\
				\noalign{\smallskip}
				\hline
				\noalign{\smallskip}
				$289$  & $3.00{\rm E}-03$    &  $3.35{\rm E}-02$  \\
				$1089$  & $8.88{\rm E}-04$    &  $1.25{\rm E}-02$    \\
				$4225$  & $2.48{\rm E}-04$    &  $7.50{\rm E}-03$   \\
				$16641$  & $1.11{\rm E}-04$    &  $2.49{\rm E}-03$    \\
				$66049$  & $8.64{\rm E}-06$    &  $1.61{\rm E}-04$     \\
				\hline
			\end{tabular}
		\end{center} 
		\caption{RMSEs and MAEs computed on Halton points and obtained by applying the BLOOCV-PU method with the IMQ as local RBF interpolant and optimal $\varepsilon_j$ (but $\delta_j=\delta$) for $f_1$.}
		\label{tabe_1a}
	\end{table}

	\subsection{Tests with non-conformal data points}
	Let us now turn into the more complex case of irregular data, which gives the results shown in Tables \ref{tabe_3a} and \ref{tabe_4a}. We use as data distribution the one shown in Figure \ref{concentric_circles_400_cerchio} (right). 
	In order to have a better understanding of their distribution we report in Table \ref{fs_dist} the two indicators of data regularity, i.e. the separation distance \eqref{sd} and  the fill distance \eqref{fd}, of the non-conformal points and we compare them with Halton data. 
	The last column gives an idea  of the quasi-uniformity constant, see Remark \ref{fill}. As evident, points coming from a Schwarz-Christoffel transformation are really far from being quasi-uniform. In this case, the interpolation process results particularly hard. In fact, we remark that the ill-conditioning  primarily grows due to the decrease of the separation distance.

	\begin{table}[ht!]
		\begin{center}
			\begin{tabular}{ccccc} 	\hline\noalign{\smallskip}
				$N$ &  data set	 & $h_{ {\cal X}_N, \Omega}$ & $q_{ {\cal X}_N}$ &  $h_{ {\cal X}_N, \Omega}$ / $q_{ {\cal X}_N}$ \\
				\hline 
				\rule[0mm]{0mm}{3ex}
				$\hskip-2pt 289$ & Halton &  $7.46{\rm E}-02$    &  $1.03{\rm E}-02$ & 	$7.20{\rm E}+00$   \\
				&  non-conformal          &  $1.87{\rm E}-01$    &  $4.90{\rm E}-03$ &  $3.81{\rm E}+01$	   \\
				\rule[0mm]{0mm}{3ex}
				$1089$     &  Halton     &  $3.93{\rm E}-02$    &  $4.33{\rm E}-03$ & 	$9.07{\rm E}+00$   \\	   
				&   non-conformal  & $1.39{\rm E}-01$    &  $1.31{\rm E}-03$ &	$1.06{\rm E}+02$  \\
				\rule[0mm]{0mm}{3ex}
				$4225$     & Halton    &  $2.19{\rm E}-02$    &  $2.19{\rm E}-03$ & 	$9.97{\rm E}+00$   \\
				&   non-conformal   & $9.94{\rm E}-02$    &  $6.84{\rm E}-04$  &	$1.45{\rm E}+02$        \\
				\hline 
			\end{tabular}
		\end{center}
		\caption{Fill and separation distances of different sets of Halton data and non-conformal points.}
		\label{fs_dist}
	\end{table}

	Since points are not well-distributed and ill-conditioning is expected, we choose as local approximant a CSRBF, specifically the Wendland $C^6$, see \eqref{WEN}.
	
	The classical PU approach is applied by fixing the shape parameter $\varepsilon=0.5$ and the radius as in \eqref{fisso}.

	\begin{table}[ht!]
		\begin{center}
			\begin{tabular}{cccc} 	\hline\noalign{\smallskip}
				$N$ & method	 & RMSE   & MAE \\
				\hline 
				\rule[0mm]{0mm}{3ex}
				$\hskip-2pt 289$ & PU & $3.28{\rm E}-02$    &  $1.90{\rm E}-01$    \\
				&  BLOOCV-PU        & $3.64{\rm E}-03$    &  $4.15{\rm E}-02$    	   \\
				\rule[0mm]{0mm}{3ex}
				$1089$     &  PU    & $1.12{\rm E}-02$    &  $2.01{\rm E}-01$  \\
				&   BLOOCV-PU   & $5.40{\rm E}-04$    &  $9.11{\rm E}-03$    \\
				\rule[0mm]{0mm}{3ex}
				$4225$     & PU   & $1.44{\rm E}-02$    &  $2.23{\rm E}-01$ 	   \\
				&   BLOOCV-PU    & $1.24{\rm E}-04$    &  $3.34{\rm E}-03$     \\
				\rule[0mm]{0mm}{3ex}
				$16641$     &  PU    & $1.12{\rm E}-02$    &  $1.92{\rm E}-01$  \\
				&   BLOOCV-PU   & $3.21{\rm E}-05$    &  $7.05{\rm E}-04$         \\
				\rule[0mm]{0mm}{3ex}
				$66049$     &  PU   & $1.26{\rm E}-02$    &  $2.45{\rm E}-01$   \\
				&  BLOOCV-PU    & $1.14{\rm E}-05$    &  $3.70{\rm E}-04$   \\
				\hline 
			\end{tabular}
		\end{center}
		\caption{RMSEs and MAEs computed on non-conformal points and obtained by using the Wendland $C^6$ as local RBF interpolant for $f_1$.}
		\label{tabe_3a}
	\end{table}    
	
	\begin{table}[ht!]
		\begin{center}
			\begin{tabular}{cccc} 	\hline\noalign{\smallskip}
				$N$ & method	 & RMSE   & MAE \\
				\hline 
				\rule[0mm]{0mm}{3ex}
				$\hskip-2pt 289$ & PU  & $5.30{\rm E}-02$    &  $5.00{\rm E}-01$  \\
				&  BLOOCV-PU        & $3.47{\rm E}-02$    &  $3.13{\rm E}-01$ 	   \\
				\rule[0mm]{0mm}{3ex}
				$1089$     &  PU    & $3.90{\rm E}-02$    &  $5.00{\rm E}-01$  \\
				&   BLOOCV-PU   & $7.11{\rm E}-03$    &  $8.38{\rm E}-02$    \\
				\rule[0mm]{0mm}{3ex}
				$4225$     & PU    & $4.63{\rm E}-02$    &  $4.99{\rm E}-01$ 	   \\
				&   BLOOCV-PU     & $2.39{\rm E}-03$    &  $4.77{\rm E}-02$      \\
				\rule[0mm]{0mm}{3ex}
				$16641$     &  PU    & $4.34{\rm E}-02$    &  $5.00{\rm E}-01$   \\
				&   BLOOCV-PU   & $7.51{\rm E}-04$    &  $8.10{\rm E}-03$      \\
				\rule[0mm]{0mm}{3ex}
				$66049$     &  PU   & $4.03{\rm E}-02$    &  $5.00{\rm E}-01$   \\
				&  BLOOCV-PU    & $8.28{\rm E}-05$    &  $1.27{\rm E}-03$   \\
				\hline 
			\end{tabular}
		\end{center}
		\caption{RMSEs and MAEs computed on non-conformal points and obtained by using the Wendland $C^6$ as local RBF interpolant for $f_2$.}
		\label{tabe_4a}
	\end{table}     
	
	In this case the classical PU method with a fixed size of the subdomains  gives inaccurate approximations. We verify with numerical experiments that this does not depend on the shape parameter; in fact, neither the use the optimal $\varepsilon$ changes the order of the approximation errors.
	The BLOOCV-PU reveals its robustness in both cases. It follows that such method with the related flexibility of the patches turns out to be meaningful when irregular data sets are considered. 
	
	\section{Applications to real world data sets}
	\label{foresta}
	This section is devoted to  test the BLOOCV-PU scheme with  two real world data sets. The first one we consider is the so-called \emph{glacier}  data set. It consists of $8345$ points  representing  digitized height contours of a glacier \cite{franke_ds,Schaback00a,Wendland01}. The difference between the highest and the lowest point  is $800$ m.
	
	As second example, we consider the so-called \emph{black forest} data set \cite{Davydov1,Davydov2}. It consists of $15885$ points representing a terrain in the neighborhood of Freiburg, Germany. In this case, the difference between the maximal and minimal heights is $1214$ m.
	
\begin{figure}
		\begin{center}
		\makebox[\textwidth]{
		\includegraphics[height=.28\textheight]{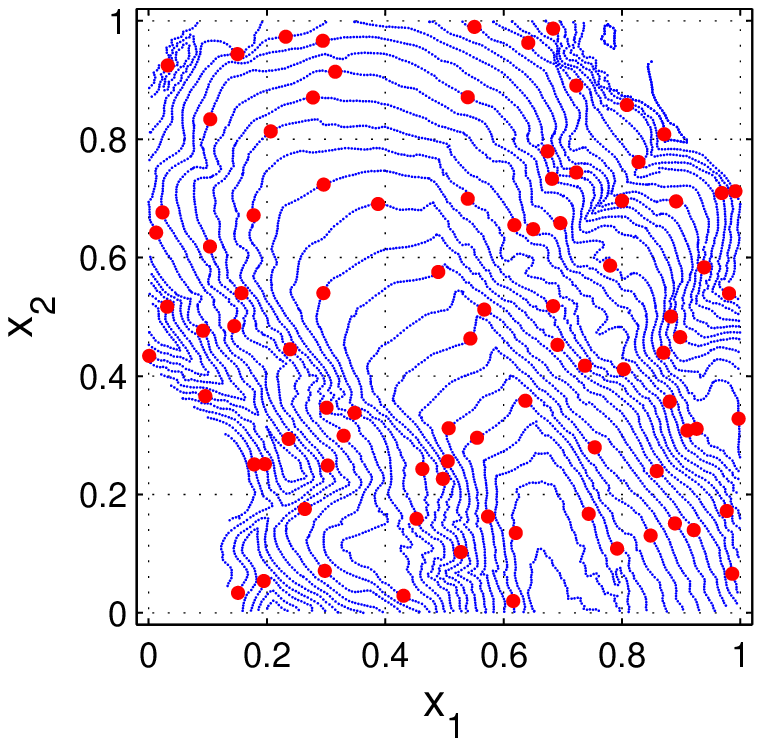}
		\hskip -1cm
		\includegraphics[height=.28\textheight]{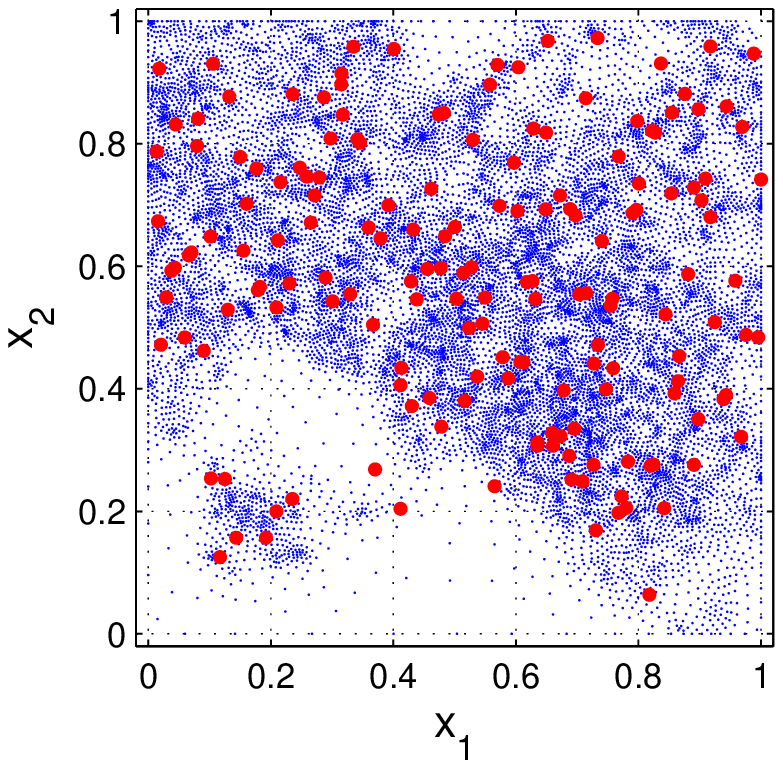}}\\
		\makebox[\textwidth]{
		\includegraphics[height=.28\textheight]{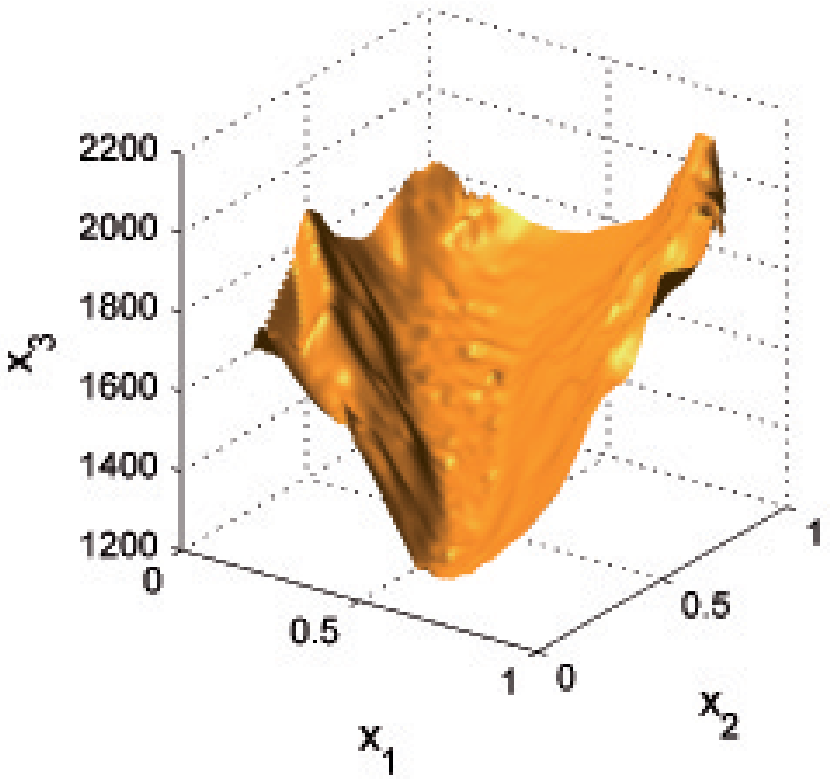}
		\hskip -1.5cm
		\includegraphics[height=.28\textheight]{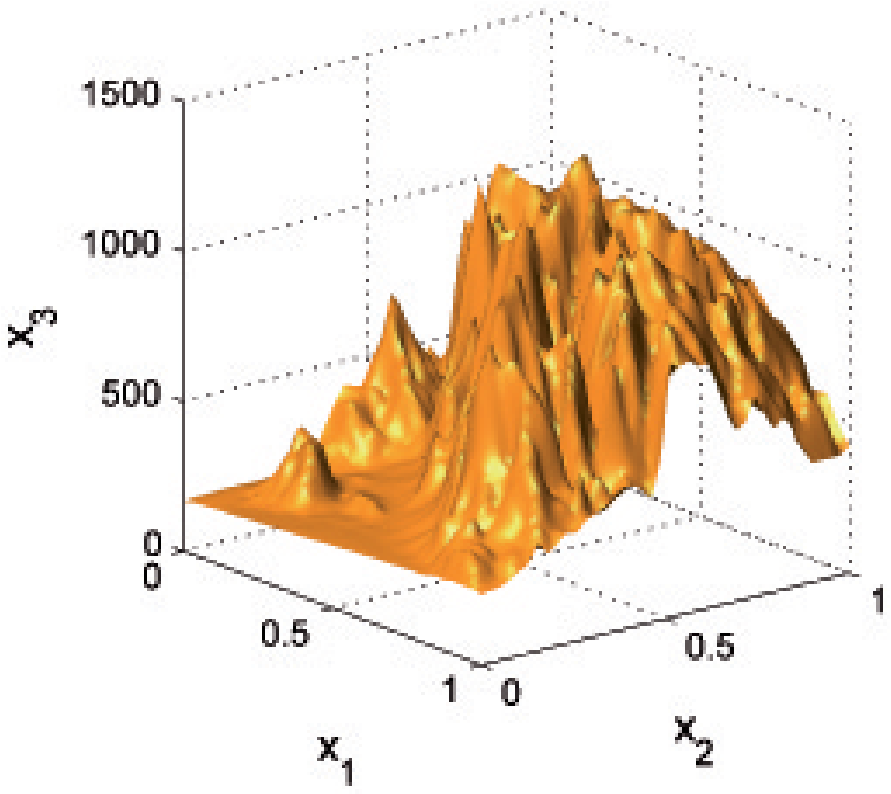}}
		\caption{A 2D view of glacier and black forest data sets (top, left to right). Graphical approximation of  glacier and  black forest data sets (bottom, left to right). The blue dots represent the set of scattered data and the red dots the points used for the validation.}
		\label{fig_dataset}
	\end{center}
	\end{figure}

	A 2D view of the data sets is plotted in Figure \ref{fig_dataset}	(top left to right, respectively).
	Because of the high variability of the points, in both cases we use  as local approximant a RBF characterized by a finite regularity, such as the Mat\'ern $C^2$, see \eqref{matern}.
	
	Figure  \ref{fig_dataset} (bottom) shows the reconstruction of the surfaces defined by the two data sets. It has been obtained evaluating the BLOOCV-PU interpolant on a grid of $80 \times 80$ points.
	
	To test the accuracy we use, as validation set,  $90$ and $170$ points of the glacier and black forest data sets, respectively. Such  data are plotted in red in Figure \ref{fig_dataset}. The errors obtained in these cases are shown in Table \ref{tabe_me}. 

	\begin{table}  
		\begin{center}
			\begin{tabular}{ccc}
				\hline\noalign{\smallskip}
				Data set  & RMSE  &  MAE  \\
				\noalign{\smallskip}
				\hline
				\noalign{\smallskip}
				Glacier & $0.65$ m    &  $3.31$ m \\
				Forest  & $5.73$ m    &  $26.0$ m\\ 
				\hline
			\end{tabular}
		\end{center} 
		\caption{RMSEs and MAEs computed on the glacier and black forest data sets and obtained by using  the BLOOCV-PU method with the  Mat\'ern $C^2$ as local RBF interpolant.}
		\label{tabe_me}
	\end{table}
	
	\section{Concluding remarks}
	\label{rem_comcl}
	
	In this paper we provided a robust tool enabling us to safely select, for each PU subdomain, both its size and the shape parameter. Numerical evidence and applications with real world measurements show that the proposed method accurately fits data with highly varying densities. Moreover, the BLOOCV-PU implementation is carried out with a new multidimensional searching procedure which has been proved to be extremely fast. 
	
	Work in progress consists in varying the shape of the PU subdomains, which here are supposed to be hyperspherical patches. This is not trivial since several requirements for the covering might be not easily satisfied.

\end{document}